\documentclass[a4paper,10pt]{article}
\usepackage[utf8x]{inputenc}
\usepackage{epsfig}
\usepackage{amssymb}
\usepackage{amsmath}
\usepackage{amsfonts}
\usepackage{moreverb}
\usepackage{graphicx}
\usepackage{rotating}
\usepackage{multirow}

\usepackage[colorlinks,bookmarksopen,bookmarksnumbered,citecolor=red,urlcolor=red]{hyperref}

\newtheorem{e-proposition}[theorem]{Proposition}

\newtheorem{e-definition}[theorem]{Definition\rm}

\title{A limitation of\\ the hydrostatic reconstruction technique\\ for Shallow Water equations}
\author{{O. Delestre$^{\dag,}$}\footnote{Laboratoire de Math\'ematiques J.A. Dieudonn\'e UMR 7351 CNRS
 UNSA \& Polytech Nice -- Sophia, Universit\'e de Nice -- Sophia Antipolis, Parc Valrose, F-06108 Nice cedex 02, France ;
 delestre@unice.fr},
 {S. Cordier}\footnote{MAPMO UMR CNRS 7349, Universit\'e Orl\'eans,
B\^atiment de math\'ematiques,
B.P. 6759, 45067 Orl\'eans cedex 2, France}, {F. Darboux}\footnote{INRA, UR 0272 Science du sol, Centre de recherche
 d'Orl\'eans, CS 40001 Ardon, F-45075 Orl\'eans cedex 2, France}
 \& {F. James$^{\dag}$}}

\begin{document}

\maketitle

\begin{abstract} English version:
 Because of their capability to preserve steady-states, well-balanced schemes for Shallow Water equations are
 becoming popular. Among them, the hydrostatic reconstruction proposed in \cite{Audusse04b}, coupled with a positive numerical
flux, allows to verify important mathematical and physical properties like the positivity of the water height and, thus, to avoid
unstabilities when dealing with dry zones. In this note, we prove that this method
exhibits an abnormal behavior for some combinations of slope, mesh size and water height.

Version Fran\c{c}aise : {\bf Une limitation de la reconstruction hydrostatique pour la r\'esolution du syst\`eme de Saint-Venant.}
De par leur capacit\'e \`a pr\'eserver les \'etats d'\'equilibre, les sch\'emas \'equilibres 
connaissent actuellement un fort d\'eveloppement dans la r\'esolution des \'equations de Saint-Venant.
En particulier, la reconstruction hydrostatique propos\'ee dans \cite{Audusse04b}, coupl\'ee \`a un flux num\'erique positif,
permet de garantir certaines propri\'et\'es comme la positivit\'e de
la hauteur d'eau et, donc, d'\'eviter certaines instabilit\'es pour traiter les zones s\`eches.
Dans cette note, nous montrons que cette m\'ethode pr\'esente un
d\'efaut pour certaines combinaisons de pente, taille de maillage et hauteur d'eau.
 \end{abstract}
%
%
%
\section*{Version fran\c{c}aise abr\'eg\'ee}
Les \'equations de Saint-Venant \eqref{eq:SW} posent des difficult\'es num\'eriques sp\'ecifiques~: pr\'eservation des \'equilibres
 stationnaires (flaques d'eau, lacs) et de la positivit\'e de la hauteur d'eau. La reconstruction hydrostatique, introduite dans
 \cite{Audusse04b,Bouchut04}, s'est impos\'ee comme une m\'ethode particuli\`erement adapt\'ee. Elle fait partie de la classe
 des sch\'emas dits \'equilibr\'es (well-balanced). Partant d'une m\'ethode de volumes finis \eqref{eq:homogeneous}, dont le flux
 num\'erique $\bf F$ est adapt\'e au syst\`eme sans topographie, on reconstruit les variables $u$, $h$, $h+z$ afin de pr\'eserver
 les \'equilibres. 
Cette m\'ethode est appliqu\'ee le cas \'ech\'eant \`a des variables d\'ej\`a reconstruites afin d'augmenter l'ordre de convergence.
La mise en \oe{u}vre compl\`ete de ce sch\'ema dans le cas d'une reconstruction lin\'eaire de type MUSCL est donn\'ee par
\eqref{hydrostatic-reconstruction}--\eqref{centered-term}.

Une reformulation de la reconstruction hydrostatique \eqref{hydrostatic-reconstruction2} met en \'evidence
 un comportement anormal pour certaines combinaisons de pente, maillage et hauteur d'eau. 
Plus pr\'ecis\'ement, il est montr\'e que le sch\'ema surestime la hauteur d'eau dans les r\'egions o\`u la relation
\eqref{critere} est v\'erifi\'ee.
 
Ce d\'efaut est illustr\'e par une s\'erie de tests num\'eriques sur une solution analytique (voir \ref{test}). Il s'av\`ere
particuli\`erement spectaculaire \`a l'ordre 1, o\`u la m\'ethode calcule une m\^eme pente apparente pour diff\'erentes valeurs
effectives (Figure~\ref{figSimulation-hydro}(Left)), mais il reste observable \`a l'ordre 2 (Figure \ref{figSimulation-hydro}(Right)).

\section{Introduction}
\label{intro}

Derived from Navier-Stokes equations, Shallow Water equations describe the water flow properties as follows
\begin{equation}
\partial_t h+\partial_x \left(hu\right)=0, \qquad
  \partial_t \left(hu\right)+\partial_x \left(hu^2+gh^2/2\right)=-gh\partial_x z,\label{eq:SW}
\end{equation}
where the unknowns are the water height $h(t,x)$, the velocity $u(t,x)$, and the topography $z(x)$ is a given function. In what
 follows, we note the discharge $q=hu$, the vector of conservative variables $U=(h\quad hu)^t$ and the flux
 $F(U)=(hu\quad gh^2/2+hu^2)^t$. The steady state of a lake at rest, or a puddle, ($h+z=Cst$ and $u=q=0$) is a particular solution
 to \eqref{eq:SW}. Since \cite{Bermudez94}, it is well known that the topography source term needs a special treatment in
 order to preserve at least this equilibrium. Such schemes are said to be well-balanced (since \cite{Greenberg96}).

In the following, we present briefly the so-called hydrostatic reconstruction method, which permits,
when coupled to a positive numerical flux, to obtain a family of well-balanced schemes that can
 preserve the water height nonnegativity and deal with dry zones.
We show that this method, presented in \cite{Audusse04b,Bouchut04} and widely used, 
 fails for some combinations of slope, mesh size and water height.
We give the criteria that ensures the accuracy of the results.

 \section{The numerical method} \label{method}

The hydrostatic reconstruction follows the general principle of reconstruction methods. We start from a first order finite volume
 scheme for the homogeneous form of system \eqref{eq:SW}: choosing a numerical flux ${\bf F}(U_L,U_R)$ ({\it e.g.}
Rusanov, HLL, VFRoe-ncv, kinetic), a finite volume scheme takes the general form
	\begin{equation}
U_i^*=U_i^n-\dfrac{\Delta t}{\Delta x} \left[{\bf F}(U_i,U_{i+1})-{\bf F}(U_{i-1},U_i)\right],\label{eq:homogeneous}
	\end{equation}
where $\Delta t$ is the time step and $\Delta x$ the space step. Now the idea is to modify this scheme by applying the flux to
 {\sl reconstructed variables}. Reconstruction can be used to get higher order schemes, in that case higher order in time is achieved
 through TVD-Runge-Kutta methods. 
The aim of the hydrostatic reconstruction, which is described in the next section, is to be well-balanced, in the sense
 that it is designed to preserve at least steady states at rest ($u=0$). When directly applied on the initial scheme, it leads to
order one methods, while coupling it with high order accuracy reconstruction increases the order.

\subsection{The hydrostatic reconstruction}\label{hydro}

We describe now the implementation of this method for high order accuracy. The first step is to perform
a high order reconstruction (MUSCL, ENO, UNO, WENO). 
 To deal properly with the topography source term $\partial_x z$, this reconstruction has to be
 performed on \(u\), \(h\) and \(h+z\), see \cite{Bouchut04}. This gives us the
 reconstructed values \((U_-,z_-)\) and \((U_+,z_+)\), on which we apply the hydrostatic reconstruction
 \cite{Audusse04b,Bouchut04} on the water height, namely
	\begin{equation}
 \left\{\begin{array}{l}
         h_{i+1/2L}=\max(h_{i+1/2-}+z_{i+1/2-}-\max(z_{i+1/2-},z_{i+1/2+}),0),\\
         U_{i+1/2L}=(h_{i+1/2L},h_{i+1/2L}u_{i+1/2-}),\\
         h_{i+1/2R}=\max(h_{i+1/2+}+z_{i+1/2+}-\max(z_{i+1/2-},z_{i+1/2+}),0),\\
         U_{i+1/2R}=(h_{i+1/2R},h_{i+1/2R}u_{i+1/2+}).
        \end{array}\right.\label{hydrostatic-reconstruction}
	\end{equation}

 With the hydrostatic reconstruction, the finite volume scheme \eqref{eq:homogeneous} is modified as follows
\begin{equation}
 U^*_i=
U^n_i-\Delta t \Phi(U^n)=
U^n_i-\dfrac{\Delta t}{\Delta x}\left[F^n_{i+1/2L}-F^n_{i-1/2R}-Fc^n_i\right],\label{discretization}
\end{equation}
where
\begin{equation}\label{flux}
 F^n_{i+1/2L}=F^n_{i+1/2}+S^n_{i+1/2L},\qquad
F^n_{i-1/2R}=F^n_{i-1/2}+S^n_{i-1/2R}
\end{equation}
are left (respectively right) modifications of the initial numerical flux for the homogeneous problem. In this formula
the flux is now applied to reconstructed variables:
$F^n_{i+1/2}={\bf F}(U^n_{i+1/2L},U^n_{i+1/2R})$, and we have introduced
\begin{equation}
 S_{i+1/2L}^n=\left(\begin{array}{c}
                0\\
	      \dfrac{g}{2}(h^2_{i+1/2-}-h^2_{i+1/2L})
               \end{array}\right),\quad
 S_{i-1/2R}^n=\left(\begin{array}{c}
                0\\     
             \dfrac{g}{2}(h^2_{i-1/2+}-h^2_{i-1/2R})
       \end{array}\right).\label{flux-modification}
\end{equation}
Finally, a centered source term has to be added to preserve consistency and well-balancing (see \cite{Audusse04b,Bouchut04}):
\begin{equation}
 Fc_i=\left(\begin{array}{c}
             0\\
	    -g\dfrac{h_{i-1/2+}+h_{i+1/2-}}{2}(z_{i+1/2-}-z_{i-1/2+})
            \end{array}\right).
\label{centered-term}
\end{equation}
Formula \eqref{centered-term} preserves the second order accuracy, but has to be modified for higher order approximations.

\subsection{The hydrostatic reconstruction rewritten}

\begin{figure}
\begin{tabular}{ccc}
\includegraphics[width = 0.47 \textwidth]{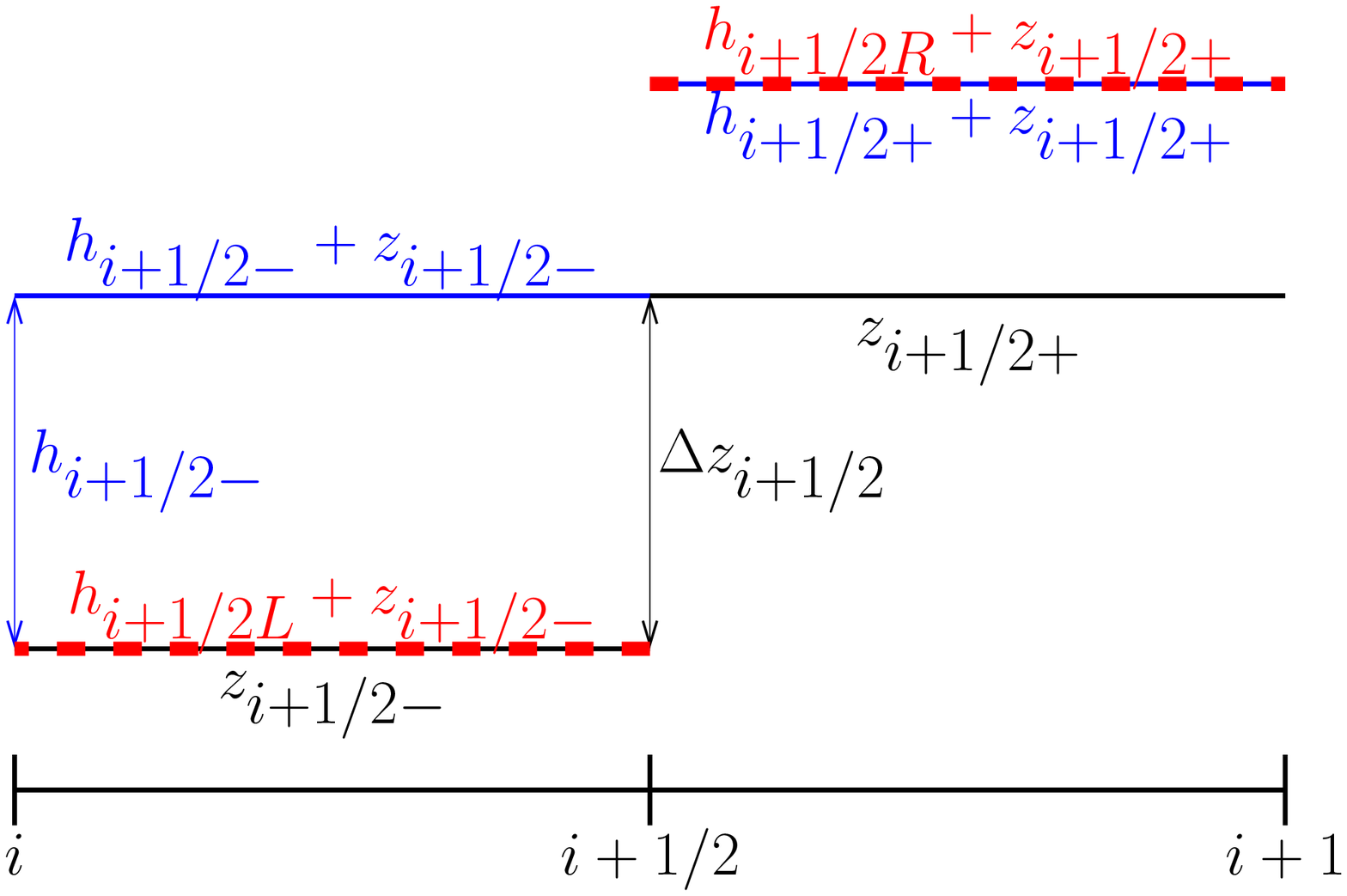}&\hspace*{2em}&
\includegraphics[width = 0.47 \textwidth]{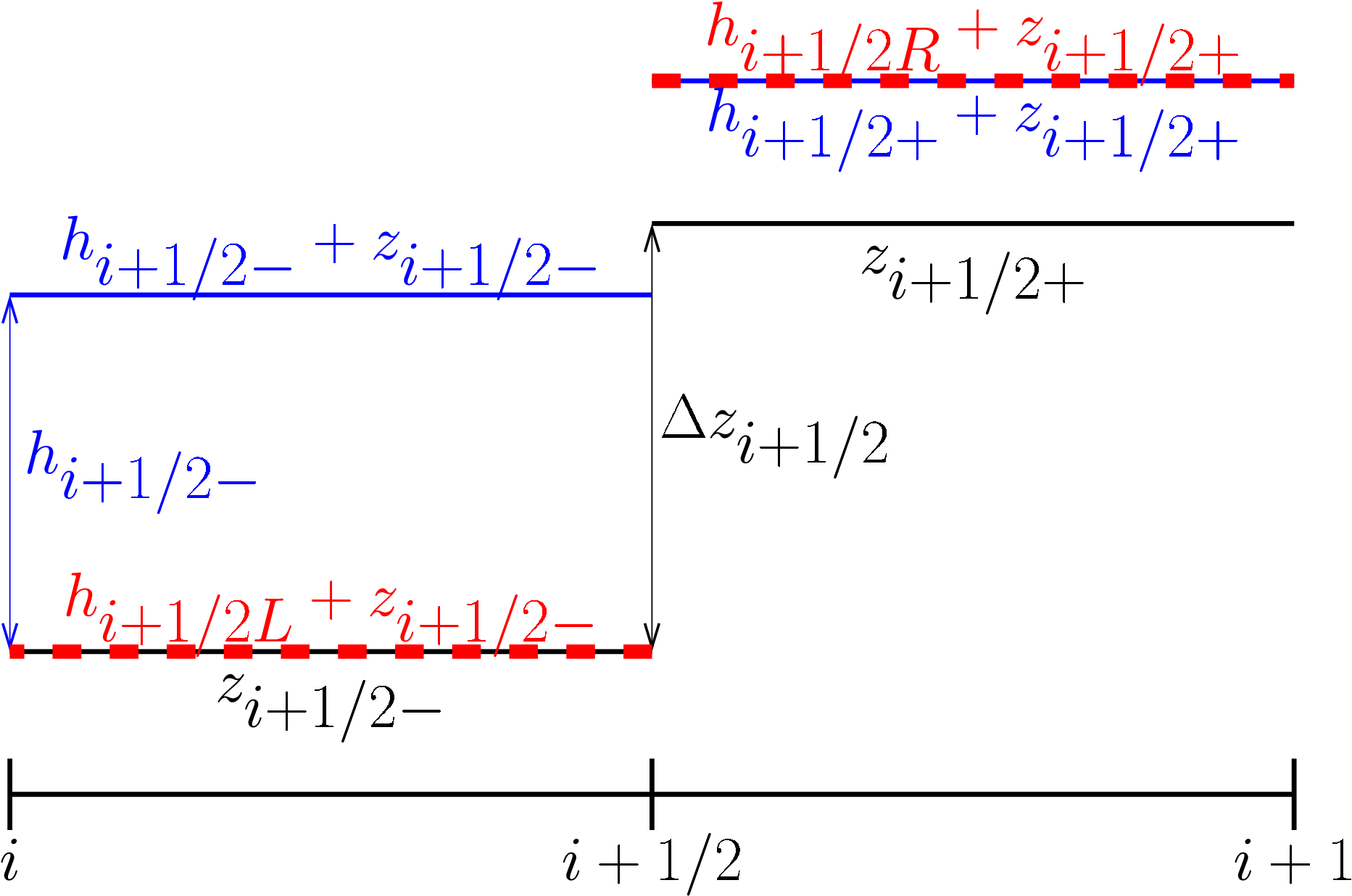}
\end{tabular}
\caption{Default of the hydrostatic reconstruction. Left: threshold non activated (limit case). Right: threshold activated}
\label{figIllustration-hydro}
\end{figure}

We define $\Delta z_{i+1/2}= z_{i+1/2+}-z_{i+1/2-}$. Once the space step and the high order reconstruction chosen, this is a
fixed sequence. Now, the reconstructed variables \eqref{hydrostatic-reconstruction} write
	\begin{equation}
\left\{\begin{array}{l}
         h_{i+1/2L}=\max(h_{i+1/2-}-\max(0,\Delta z_{i+1/2}),0),\\
         h_{i+1/2R}=\max(h_{i+1/2+}+\min(0,\Delta z_{i+1/2}),0).
\end{array}\right.\label{hydrostatic-reconstruction2}
	\end{equation}
With \eqref{hydrostatic-reconstruction2}, the defect of the hydrostatic reconstruction becomes apparent. To fix the ideas,
suppose the local slope is positive, hence $\Delta z_{i+1/2}>0$, as in Figure \ref{figIllustration-hydro}.
Then, for all $\Delta z_{i+1/2}$ such that $\Delta z_{i+1/2}\geq h_{i+1/2-}\geq 0$, 
the reconstruction gives $h_{i+1/2L}=0$, while $h_{i+1/2R}$ remains unchanged. In that case, the reconstruction prevents
an unphysical negative value for $h_{i+1/2L}$, the counterpart of this being an underestimated difference $h_{i+1/2L}-h_{i+1/2R}$. 
Therefore, there is a lack in the numerical flux computed from the modified Riemann problem, which gives an underestimated
velocity and consequently an overestimated height. This can be interpreted in terms of reconstructed slope as well, which
is underestimated. 

In the general case, we can write the following local criterion of ``non validity'' for the hydrostatic reconstruction.
\begin{e-proposition}
For a fixed discretization, if for some $i_0\le i \le i_1$ one has
	\begin{equation}\label{critere}
\Delta z_{i+1/2}\geq h_{i+1/2-}\geq 0, \qquad\mbox{or}\quad -\Delta z_{i+1/2}\geq h_{i+1/2+}\geq 0,
	\end{equation}
then the hydrostatic reconstruction will overestimate (resp. underestimate) the water height (resp. velocity).
\end{e-proposition}

Notice that from a theoretical viewpoint, this is not limiting. Indeed, since this class of schemes is consistent with the 
system of partial differential equations (see \cite{Bouchut04}), the problem disappears when refining the discretization.
However, it has to be taken into account for practical computations,  with a fixed discretization. It is particularly
apparent for order $1$ schemes, but also remains present for order $2$.

\section{Numerical illustration}\label{test}

For numerical illustration purpose, we introduce an explicit solution which
consists in a supercritical steady flow on an inclined plane with constant slope $\partial_x z=\alpha$
(it is referenced in the SWASHES library \cite{Delestre}). Steady states are solutions to
\begin{equation*}
        q(x,t)=q_0=Cst,\qquad
	\partial_x \left(hu^2+g\dfrac{h^2}{2}\right)=-gh\partial_x z,
\end{equation*}
and the height profile $h(x)$ must be a solution to Bernoulli's law rewritten as a third order equation in $h$:
	\begin{equation}
 h^3+h^2\left(\alpha x- \dfrac{q_0^2}{2gh_0^2}-h_0\right)+\dfrac{q_0^2}{2g}=0,\label{polynom}
	\end{equation}
where $h_0$ and $q_0=h_0u_0$ are the height and discharge at $x=0$, which completely determine the profile 
since the flow is supercritical. A careful study of the roots of the polynomial shows that the supercritical height profile
is decreasing in $x$: $h(x)\le h_0$ for all $x\ge 0$.

The numerical strategy we choose consists in the HLL flux and a modified MUSCL reconstruction. In \cite{Delestre10}, this
 combination of flux and linear reconstruction has shown to be the best compromise between accuracy, stability and CPU
 time cost. We refer to \cite{Bouchut04,Delestre10} for a presentation of the HLL flux.
The MUSCL reconstruction of a scalar function \(s\in\mathbb{R}\) is
\begin{equation}
 s_{i-1/2+}=s_i-\dfrac{\Delta x}{2}Ds_i, \qquad 
s_{i+1/2-}=s_i+\dfrac{\Delta x}{2}Ds_i,\label{linear-reconstruction}
\end{equation}
where the ``minmod'' operator $D$ is given by
\begin{equation}\label{modif}
 Ds_i=\text{minmod}\left(\dfrac{s_i-s_{i-1}}{\Delta x},\dfrac{s_{i+1}-s_i}{\Delta x}\right), \quad
 \text{minmod}(x,y)=\begin{cases}
                      \min(x,y)&  \text{if}\; x,y\geq 0, \\
		      \max(x,y)&  \text{if}\; x,y\leq 0, \\
		      0  & \text{else.}
                     \end{cases}
\end{equation}
In order to keep the discharge conservation, the reconstruction of the velocity has to be modified as
\begin{equation*}
 u_{i-1/2+}=u_i-\dfrac{h_{i+1/2-}}{h_i}\dfrac{\Delta x}{2}Du_i\qquad
u_{i+1/2-}=u_i+\dfrac{h_{i-1/2+}}{h_i}\dfrac{\Delta x}{2}Du_i.
\end{equation*}
Notice that if we take $Ds_i=0$ in \eqref{linear-reconstruction}, then $z_{i+1/2-}=z_{i+1/2-}=z_i$ so that the centered term
 \eqref{centered-term} disappears, and we recover the first order scheme in space.
Second order in time is achieved through a classical Heun predictor-corrector method.

We turn now to the specific data for simulations. The domain length is $L=10\;\text{m}$, and we choose initial data
	\begin{equation*}
h_0=h(x=0,t)=0.02 \;\text{m},\qquad q_0=q(x=0,t)=0.01\;\text{m}^2/\text{s},
	\end{equation*}
which is indeed supercritical.
All simulations are performed with a space step $\Delta x=0.1\;\text{m}$, the time step is fixed in order to satisfy
 the CFL condition.

\begin{figure}
\begin{tabular}{ccc}
\includegraphics[angle=-90,width = 0.47 \textwidth]{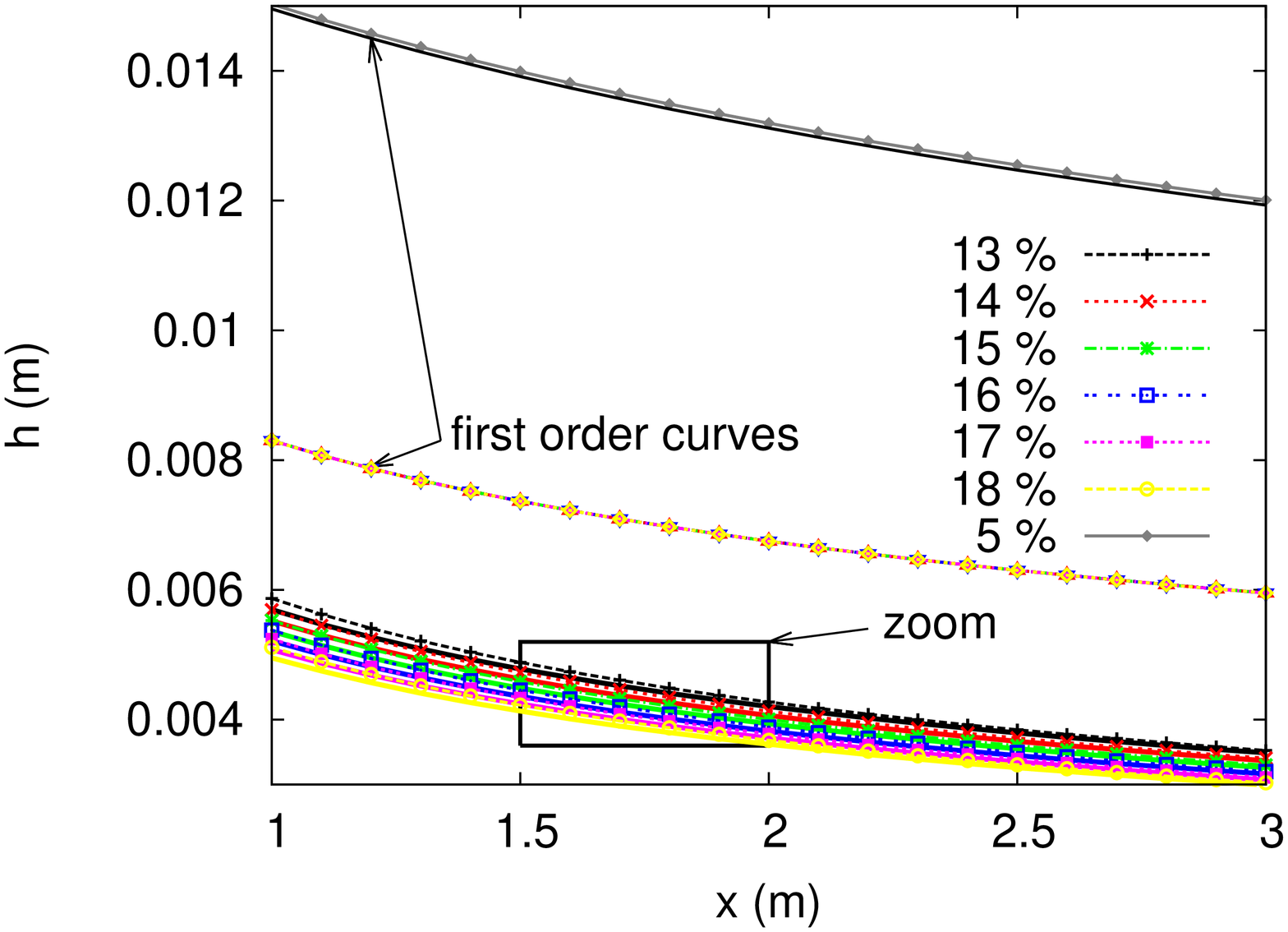} & \hspace*{2em}&
\includegraphics[angle=-90,width = 0.47 \textwidth]{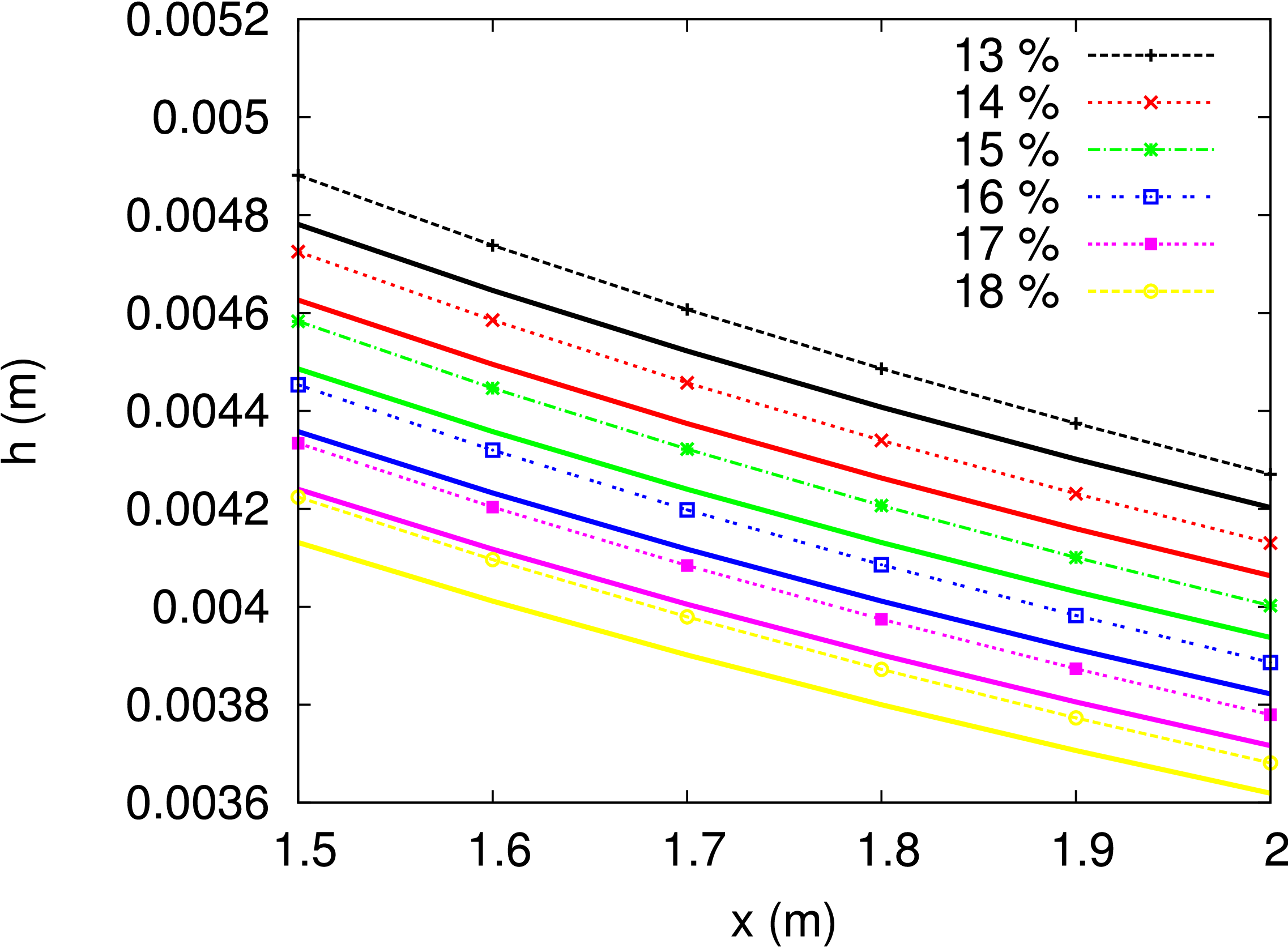}
\end{tabular}
\caption{Default of hydrostatic reconstruction: water height at first and second order of accuracy for
different slopes. Dotted curves are simulations, plain curves are exact solutions.
Left:~first and second order curves.
Right:~zoom on second order curves}
\label{figSimulation-hydro}
\end{figure}

The analytical solution is computed with $7$ negative slopes $\alpha=5\%$, $13\%$, $14\%$, $15\%$, $16\%$, $17\%$ and $18\%$,
both at first and second order.
Numerical results are compared to the analytical results in Figure~\ref{figSimulation-hydro}, where a part of the domain
is displayed ($x$ between $1$ m and $3$ m).

With this set of data, a domain of admissible slopes can be estimated for the order $1$ hydrostatic reconstruction. Indeed since
 $\Delta z_{j+1/2} = |\alpha| \Delta x$, inserting a characteristic height of the flow $h^*$ in \eqref{critere} gives a bound for
 the slope.
One can use $h^*=h_0$, the incoming height. Since the height profile is decreasing, 
the whole profile will be wrong if $h_0\le|\alpha|\Delta x$, that is here $|\alpha|\ge 20\%$. But actually, since the height decreases
quite rapidly, a more accurate estimate is obtained for $x=1.5$ m with $h^*=0.005$ m (see Fig. \ref{figSimulation-hydro}), which
 leads to slopes above $5\%$.

We observe on Figure~\ref{figSimulation-hydro}(Left) that with the first order scheme, the effect of the hydrostatic reconstruction
is so important that all curves for slopes between $13\%$ and $18\%$ are superposed. In that case, the apparent result
is the simulation of a single slope, namely $13\%$. For a $5\%$ slope, we still observe a slightly overestimated water height, as
 anticipated since $5\%$ is a limit case as observed above.
With the second order, the water heights are still overestimated, but in a very slighter way, and the different curves are no
longer identical (Figure~\ref{figSimulation-hydro}(Right)).

\section{Conclusions}\label{conclusions}
The hydrostatic reconstruction may fail for certain combinations of water height, slope and mesh size,
namely in regions where \eqref{critere} holds. The defaults are particularly apparent for order $1$ schemes, leading
to wrongly estimated slopes, but still remain at order $2$, with some overestimated water heights. 
We emphasize that the problem disappears when refining the mesh,
but has to be taken into account for a given discretization. 
The generalization of the hydrostatic reconstruction proposed in \cite{castro} does not exhibit the limitation
 discussed here, even for first order schemes, but positivity is not ensured. 
Other schemes involving threshold values ({\it e.g.} \cite{Bouchut10,Liang09}) very likely encounter the same kind of problem.
Alternatively, the scheme recently introduced in \cite{Berthon} preserves the water height positivity and does not suffer 
from this problem. Notice finally that criterion \eqref{critere} may be of some utility for adaptive mesh schemes, such as
the ones used in Gerris \cite{popinet}.

\section*{Acknowledgments}
This study is part of the ANR project METHODE \#ANR-07-BLAN-0232
granted by the French National Agency for Research.
The authors wishes to thank M.-O. Bristeau, Ch. Berthon and M.-H. Le for fruitful discussions.

\end{document}